\documentclass[11pt]{article}
\usepackage{epsfig}
\usepackage{amssymb}
\usepackage{bm}
\newcommand{\be}{\begin{equation}}
\newcommand{\ee}{\end{equation}}
\newcommand{\bea}{\begin{eqnarray}}
\newcommand{\eea}{\end{eqnarray}}
\newcommand{\binom}[2]{{#1 \choose #2}}
\newcommand{\mbf}[1]{\mbox{\boldmath$#1$}}

\begin{document}

%\DeclareGraphicsExtensions{.jpkg,.pdf,.png}

\title{Sylvester-Cayley vector partitions algorithm and the
Gaussian polynomials}
\author{Boris Y. Rubinstein\\
Stowers Institute for Medical Research
\\1000 50th St., Kansas City, MO 64110, U.S.A.}
\date{\today}

\maketitle
%%%%%%%%%%%%%%%%%%%%%%%%%%%%%%%%%%%%%%%%%%%%%%%%%%%%%%%%%%%%%%%%%%%%%
\begin{abstract}
We extend an %elimination
algorithm suggested 
in 1858 by Sylvester \cite{Sylv1} and 
implemented in 1860 by Cayley \cite{Cayley1860}
for a problem of double partitions and apply it to 
derivation of explicit expressions for coefficients of the
Gaussian polynomials through convolution of restricted partition functions.
\end{abstract}


%%%%%%%%%%%%%%%%%%%%%%%%%%%%%%%%%%%%%%%%%%%%%%%%%%%%%%%%%%%%%%%%%%%%%
%\newpage
%\section{Introduction}
%\label{intro}

\section{Partition functions}
%The problem of integer partition into smaller integers is equivalent to

%%%%%%%%%%%%%%%%%%%%%%%%%%%%%%%%%%%
%%%   Scalar partitions
%%%%%%%%%%%%%%%%%%%%%%%%%%%%%%%%%%%%

\subsection{Scalar partition function}
\label{ScalarPartitionFunction}
Consider a problem of 
finding a number of nonnegative integer solutions of a Diophantine equation
\be
\sum_{i=1}^m x_i d_i = {\bf x}\cdot{\bf d} = s,
\quad
{\bf d} = \{d_1,d_2,\ldots,d_m\}.
\label{coin1}
\ee
A restricted
partition function $W(s,{\bf d})$
solving the above problem is a
number of partitions of an integer $s$ into positive integers
$\{d_1,d_2,\ldots,d_m\}$,
each not greater than $s$. The generating function for $W(s,{\bf d})$
has a form
\be
G(t,{\bf d})=\prod_{i=1}^m\frac{1}{1-t^{d_{i}}}
 =\sum_{s=0}^{\infty} W(s,{\bf d})\;t^s\;.
\label{WGF}
\ee
%where $W(s,{\bf d})$ satisfies a basic recursive relation
%\be
%W(s,{\bf d}) - W(s-d_m,{\bf d}) = W(s,{\bf d}')\;,
%\quad
%{\bf d}' = \{d_1,d_2,\ldots,d_{m-1}\}.
%\label{Wrecurs}
%\ee
In other words the partition function equals a constant term in Taylor expansion
of the following expression
\be
W(s,{\bf d}) = \mbox{const}_t
\left[
t^{-s}\prod_{i=1}^m(1-t^{d_{i}})^{-1}
\right].
\label{const1}
\ee
Sylvester proved \cite{Sylv2} a statement about splitting of the partition
function into periodic and non-periodic parts and showed that the
restricted partition function may be presented as a sum of the {\em Sylvester waves}
\cite{Rub04}
\be
W(s,{\bf d}) = \sum_{j=1} W_j(s,{\bf d})\;,
\label{SylvWavesExpand}
\ee
where summation runs over all distinct factors $j$
of the elements $d_i$.
The wave $W_j(s,{\bf d})$ is a quasipolynomial in $s$
closely related to prime roots $\rho_j$ of unity;
it is a coefficient of
${t}^{-1}$ in the series expansion in ascending powers of $t$ of
the generator
\be
F_j(s,t)=
\sum_{\rho_j} \frac{\rho_j^{-s} e^{st}}{\prod_{k=1}^{m}
        \left(1-\rho_j^{d_k} e^{-d_k t}\right)}\;.
\label{generatorWj}
\ee
The summation is made over all prime roots of unity
$\rho_j=\exp(2\pi i n/j)$ for $n$ relatively prime to $j$
(including unity) and smaller than $j$.
It was shown \cite{Rub04} that it is possible to
express the Sylvester wave as a finite sum of the Bernoulli polynomials of
higher order.

\subsection{Vector partition function}
%\label{VectorPartitionFunction}
Consider a function $W({\bf s},{\bf D})$ counting the number of integer
nonnegative
solutions ${\bf x} \ge 0$
to the linear system ${\bf D} \cdot {\bf x} = {\bf s}$, where
${\bf D}$ is a nonnegative integer $l \times m$ matrix.
The {\it vector partition function} (VPF) $W({\bf s},{\bf D})$ is
a natural generalization of the restricted partition function to the
vector argument.

The generating function for the VPF reads
\be
G({\bf t},{\bf D})=\prod_{i=1}^m \frac{1}{1-{\bf t}^{{\bf c}_i}} =
\sum_{{\bf s}} W({\bf s},{\bf D}) {\bf t^s},
%\sum_{{\bf s}} W({\bf s},\{{\bf c}_1,\ldots,{\bf c}_m\}) {\bf t^s}\;,
\quad
{\bf t^s} = \prod_{k=1}^l t_k^{s_k},
\quad
{\bf t}^{{\bf c}_i} = \prod_{k=1}^l t_k^{d_{ki}},
\label{WvectGF}
\ee
where ${\bf c}_i,\ 1 \le i \le m,$
denotes the $i$-th column of the matrix ${\bf D}$.
%We also use another representation of the matrix ${\bf D}$ via
%its rows ${\bf r}_k,\ (1\le k \le l),$ so that the element $D_{ki}$
%is given by $i$-th element of ${\bf r}_k$ and $k$-th element of
%${\bf c}_i$.
Note that some elements $d_{ki}$ might equal zero.

%A method of computation of the VPF $W({\bf s},{\bf D})$
%was suggested \cite{Rub06} being
%a direct generalization of the approach developed in \cite{Rub04}.
%To this end {\em vector} Bernoulli and Eulerian polynomials
%of higher order were introduced to find explicit
%expression for $W({\bf s},{\bf D})$. One of the drawbacks of this
%method is that it does not provide any mechanism to define
%boundaries of so called {\em chambers} \cite{Sturmfels1995} --
%regions in $m$-dimensional space where the partition function
%is characterized by a specific expression.

The partition problems that include inequalities can be
reduced to the VPF problem as follows.
Assume that we have a system of $k$ linear equations and
$l-k$ linear inequalities for $m > l$
integer nonnegative variables $x_i$:
\be
\sum_{i=1}^m d_{ji} x_i = s_j,\ (1 \le j \le k),
\quad
\sum_{i=1}^m d_{ji} x_i \le n_j,\ (k+1 \le j \le l).
\label{systorig}
\ee
Introduce $l-k$ new integer variables
$x_{i},\ (k+1 \le i \le l)$ that enter additively
all above equations and inequalities with corresponding
factors $d_{ji} = 0$ for $1 \le j \le k$, and
$d_{ji} = 1$ for $k+1 \le j \le l$, and
transform all inequalities into equations
\be
\sum_{i=1}^{m+l-k} d_{ji} x_i = s_j,\ (1 \le j \le k),
\quad
\sum_{i=1}^{m+l-k} d_{ji} x_i = n_j,\ (k+1 \le j \le l).
\label{systeqonly}
\ee
It is easy to see that the number of integer
nonnegative solutions of (\ref{systeqonly}) equals
that of (\ref{systorig}), and thus any system of linear equations
with constraints can be reduced to the
standard vector partition problem.

\subsection{Sylvester-Cayley method of compound partitions}
The problem of vector partitions has a long history and J.J. Sylvester 
made a significant contribution to its solution.
He put forward \cite{Sylv1} an idea to reduce VPF
into a sum of scalar partitions (\ref{WGF}) (see also \cite{Mathews1897}).
The reduction is an iterative process based on the 
variable elimination in the generating function (\ref{WvectGF}).
Sylvester considered a specific double partition problem as 
an illustration of his method and determined regions
({\it chambers}) of a plane $\{s_1,s_2\}$ each having a unique
expression for VPF valid in this region only.
He showed that the 
expressions in the adjacent chambers coincide at their
common boundary (see also \cite{Sturmfels1995}).

This approach was successfully applied by 
A.Cayley \cite{Cayley1860} to double partitions subject
to some restrictions on the {\it positive} elements of matrix ${\bf D}$ 
(the columns ${\bf c}_i$ are linearly independent and 
the inequality $c_{i2} < s_2+2$ should hold for all $1\le i \le m$).

In this manuscript we present the result for double partitions 
obtained by Cayley and 
extend it to apply to the Gaussian polynomial problem.

%%%%%%%%%%%%%%%%%%%%%%%%%%%%%%%%%%%
%%%   Douple partitions
%%%%%%%%%%%%%%%%%%%%%%%%%%%%%%%%%%%%

\section{Double partitions}

Introduce an augmented matrix ${\bf E}$ obtained by prepending
the column vector ${\bf s}$ to the matrix ${\bf D}$,
so that $e_{10}=s_1=r$ and $e_{20}=s_2=\rho$.
Denote the rest elements of first row of ${\bf E}$ as
$e_{1i}=d_{1i}=b_i,\ (1 \le i \le m)$, and of the
second row by $e_{2i}=d_{2i}=\beta_i,\ (1 \le i \le m)$.
For simplicity assume that all $\beta_i > 0$,
and perform a partial fraction expansion (PFE) step to present
$G({\bf t},{\bf D})$ as sum of $m$ fractions
\be
G({\bf t},{\bf D}){\bf t}^{-{\bf s}} ={\bf t}^{-{\bf s}}
\prod_{i=1}^m \frac{1}{1-t_1^{b_i}t_2^{\beta_i}} =
\sum_{i=1}^m T_i({\bf t}),
\quad
T_i = \frac{A_i({\bf t}){\bf t}^{-{\bf s}}}{1-t_1^{b_i}t_2^{\beta_i}}.
\label{mCayley1}
\ee
Then the solution reads as a sum of the terms %$U_i$ where
\be
U_i = \mbox{const}_{\bf t}
\left[
T_i({\bf t})
\right].
\label{mCayley2}
\ee
Cayley showed that the first term $U_1$
corresponding to elimination of the second column of
matrix ${\bf E}$
can be written as
\be
U_1 = \mbox{const}_{t_1}
\left[
t_1^{-(r\beta_1 - b_1\rho)}
\prod_{j\ne1}^m (1-t_1^{b_j\beta_1-b_1\beta_j})^{-1}
\right].
\label{mCayley3}
\ee
Comparing this to (\ref{const1}) we see that
$U_1$ corresponds to $W(r\beta_1 - b_1\rho,{\bf d}_1)$,
where the elements $d_{1j}$ of the vector ${\bf d}_1$ are given by
$d_{1j} = b_j\beta_1-b_1\beta_j,\ j\ne 1$.
Effectively, this result corresponds to elimination \cite{Dickson2}
of the first unknown $x_1$ in the system ${\bf D}\cdot{\bf x} = {\bf s}$
that leads to a single Diophantine equation 
${\bf d}_1 \cdot {\bf x}' = r\beta_1 - b_1\rho$,
where ${\bf x}'$ is the vector of unknowns obtained 
from ${\bf x}$ by dropping $x_1$.

Similar transformations allow to find all $m$ contributions
to obtain
\be
W({\bf s},{\bf D}) = \sum_{i=1}^m U_i =
\sum_{i=1}^m W(L_i,{\bf d}_i),
\quad
L_i = r\beta_i - b_i\rho,
\quad
d_{ij} = b_j\beta_i-b_i\beta_j,\ j\ne i,
\label{mCayley8}
\ee
where both $L_i$ and $d_{ij}$ as the determinants
of $2\times2$ matrices made of the columns $\{{\bf c}_0,\ {\bf c}_i\}$ and
$\{{\bf c}_j,\ {\bf c}_i\}$, respectively.
Cayley mentioned \cite{Cayley1860} that the contribution of each term
in (\ref{mCayley8}) is nonzero only when $L_i$ is nonnegative, thus
%he was the first one who effectively 
reintroducing the notion of
{\em chambers} of the vector partition.
He also pointed out that when some elements of the vector
${\bf d}_i$ are negative (say, $d_{ij_k}<0$ for $1 \le k \le K$)
one has
\be
W(L_i,{\bf d}_i) =
(-1)^{K} W(L_i - \sum_{k=1}^{K} |d_{ij_k}|,|{\bf d}_i|),
\quad
|{\bf d}_i| = \{|d_{ij}|\}.
\label{mCayley8a}
\ee

%\subsection{Matrix with zero elements in one row}
Consider now a case when one or more elements in one 
(say, the second) row of the matrix
${\bf D}$ are zeroes. Denote the number of nonzero
elements as $m_0$, and without loss
of generality we  assume that $\beta_i \ne 0$ for
$1 \le i \le m_0$, and  $\beta_i = 0$ for
$m_0+1 \le i \le m$.
Using the above result for the columns with $\beta_i = 0$
we observe that each of them leads to a contribution
$W(L_i,{\bf d}_i)=0$ that vanishes as its argument $L_i = - b_i\rho < 0$ is negative.
Thus we obtain
\be
W({\bf s},{\bf D}) =
\sum_{i=1}^{m_0} W(L_i,{\bf d}_i),
\quad
L_i = d_{i0} = r\beta_i - b_i\rho,
\quad
d_{ij} = b_j\beta_i-b_i\beta_j,\ j\ne i.
\label{mCayley80}
\ee

%%%%%%%%%%%%%%%%%%%%%%%%%%%%%%%%%%%
%%%   Gaussian polynomial
%%%%%%%%%%%%%%%%%%%%%%%%%%%%%%%%%%%%

%%%%%%%%%%%%%%%%%%%%%%%%%%%%%%%%%%%%%%%%%%%%%%%%%
\section{Gaussian polynomials and their coefficients}
%%%%%%%%%%%%%%%%%%%%%%%%%%%%%%%%%%%%%%%%%%%%%%%%%
The Gaussian polynomial coefficients $P_m^n(s)$ are defined as a number
of nonnegative integer solutions of a linear Diophantine equation with a constrain
\be
\sum_{i=1}^m i x_i = s,
\quad
\sum_{i=1}^m x_i \le n.
\label{gauss_orig}
\ee
The coefficients $P_m^n(s)$ of the Gaussian polynomial $G_m^n(t)$ 
satisfy \cite{Sylv3}
\be
G_m^n(t) = \frac{G_{m}(t)G_{n}(t)}{G_{m+n}(t)} = 
\sum_{s=0}^{mn} P_m^n(s) t^s,
\quad
G_m(t) = \prod_{i=1}^m\frac{1}{1-t^i},
\label{gauss_GF1}
\ee
where $G_m(t)$ denotes the generating function
for the restricted partition function 
$W_m(s) = W(s,\hat{\bf d}^m = \{1,2,\ldots,m\})$.
The quasipolynomials $P_m^n(s)$ have a finite order $mn$ and the following
properties \cite{Sylv3}
\bea
&& 
P_m^n(0)= P_m^n(mn)= 1, \quad P_m^n(s)= 0, \ \ \mbox{for} \ \ s>mn,
\nonumber \\
&&
P_m^n(s)= P_n^m(s)= P_m^n(mn-s), \quad P_m^n(mn/2-s) = P_m^n(mn/2+s),
\label{Gauss_known} \\
&&
P_m^n(s) - P_n^m(s-1) \ge 0, \ \ \mbox{for} \ \ 0\le s \le mn/2.
\nonumber 
\eea
Cayley \cite{Cayley1855} considered a problem of 
computation of the Gaussian polynomial coefficients and 
found the generating functions (\ref{gauss_GF1}) and 
(\ref{gen_func_1}, see below).
%He  also  introduced an alternative generating function
%via the following construction. For integer $p$ consider 
%expansion coefficients 
%$c_p(s)$ of the generating function
%$$
%t^{pn+p(p+1)/2} \prod_{j=0}^{m} G_{p+j}(t) = 
%\sum_{s=0} c_p(s) t^s.
%$$
%Then we have 
%$$
%P_m^n(s) = \sum_{p=0}^m (-1)^p c_p(s),
%$$
%thus we arrive at
%\be
%\sum_{p=0}^m (-1)^p 
%t^{pn+p(p+1)/2} \prod_{j=0}^{m} G_{p+j}(t) = 
%\sum_{s=0}^{mn} P_m^n(s) t^s.
%\label{gen_func_Cayley}
%\ee
%Cayley gave several examples of the coefficient computation.
Some interesting results about explicit expression of $P_m^n(s)$ 
through the regular partition functions
were reported recently in \cite{Fel16}.

\subsection{Sylvester-Cayley algorithm for Gaussian polynomial coefficients}
The problem (\ref{gauss_orig}) is a particular case of 
(\ref{systorig}) for $k=1,\ l=2$, and $d_{1i} = i,\ d_{2i}=1$,
thus it can be transformed into a VPF problem 
\be
\sum_{i=0}^m i x_i = s,
\quad
\sum_{i=0}^m x_i = n,
\label{gauss_def}
\ee
and dealt with using the Sylvester-Cayley algorithm.

Write (\ref{mCayley80}) for %the vector ${\bf s}=\{s,n}$ and the matrix
$$
{\bf D} = \left(
\begin{array}{ccccc}
1 & 1 & 1 & \ldots & 1 \\
0 & 1 & 2 & \ldots & m
\end{array}
\right),
\quad
{\bf s} = \left(
\begin{array}{c}
n \\ s
\end{array}
\right) \ ,
$$
that determines a generating function
for $P_m^n(s)$ with fixed integers $m$
and $n$
\be
%{\cal G}_m(a,t) = 
\prod_{i=0}^m \frac{1}{1-a t^i} = 
\sum_{n=0}^{\infty}G_m^n(t) a^n = 
\sum_{n=0}^{\infty}\sum_{s=0}^{mn} P_m^n(s) a^n t^s.
\label{gen_func_1}
\ee

In (\ref{mCayley80}) we have $m_0 = m,\ r=n,\ \rho=s,\ b_i=1, \beta_i = i$ and obtain
$L_i = n i-s$ and $d_{ij} = \beta_i-\beta_j = i-j$. Thus the set ${\bf d}_i$ 
can be written as a union of $\hat{\bf d}^{i}$ and $-\hat{\bf d}^{m-i}$.
%, where it is assumed that $\hat{\bf d}^{i-1} = \emptyset$.
We then observe from (\ref{mCayley80})
\be
P_m^n(s) = \sum_{i=1}^m w_i(m,n;s),
\quad
w_i(m,n;s) = W(n i - s,-\hat{\bf d}^{m-i} \cup \hat{\bf d}^{i} ),
\label{gauss_1}
\ee
where the $i$-th term in the above sum contributes for $0 \le s \le n i$. 
This immediately leads to a conclusion that $P_m^n(s) = 0$ for $s > nm$ as 
all $m$ terms in (\ref{gauss_1}) vanish.
By the representation (\ref{gauss_1}) the Gaussian polynomial $P_m^n(s)$ has $m$ chambers
bounded by lines $s=n i,\ 0 \le i\le m$.
Specifically, in the $k$-th chamber bounded by the lines 
$s=n(m-k)$ and $s=n(m-k+1)$ one has to retain only $k$ terms in the sum 
with $m-k+1 \le i \le m$.
Using (\ref{mCayley8a}) we arrive at 
\be
w_i(m,n;s) = (-1)^{m-i}
W(n i - s - s_{m-i},\hat{\bf d}^{i} \cup \hat{\bf d}^{m-i}),
\quad
s_{m} = m(m+1)/2,
\label{gauss_2}
\ee
where the partition function $W$ can be computed using the algorithm
discussed in \cite{Rub04}.
It is convenient to introduce $P_m^n(r,s)$
that explicitly accounts for the summands entering $P_m^n(s)$ in
the $r$-th chamber
\be
P_m^n(r,s) = \sum_{i=r}^m w_i(m,n;s),
\quad
(r-1)n \le s \le rn.
\label{gauss_1r}
\ee

%Note that 
%$$
%w_{m-i}(m,n;s) = (-1)^{i}
%W(nm - s - i(2n+i+1)/2,\hat{\bf d}^{i} \cup \hat{\bf d}^{m-i}).
%$$

\subsection{Convolution of restricted partitions}
%It is also possible to find an alternative representation of the result
%(\ref{gauss_1},\ref{gauss_2}).
Consider a general term $W(s,\hat{\bf d}^{k_1} \cup \hat{\bf d}^{k_2})$ 
in (\ref{gauss_2}).
Its generating function reads
$G(t) = G_{k_1}(t) G_{k_1}(t)$ and we obtain
$$
\sum_{s=0} W(s,\hat{\bf d}^{k_1} \cup \hat{\bf d}^{k_2}) t^s = 
\sum_{s_1=0} \sum_{s_2=0} W_{k_1}(s_1)  W_{k_2}(s_1)t^{s_1+s_2},
$$
leading to a discrete convolution (the Cauchy product)
$$
W(s,\hat{\bf d}^{k_1} \cup \hat{\bf d}^{k_2}) = 
\sum_{k=0}^{s} W_{k_1}(k) W_{k_2}(s-k),
$$
and setting here $k_1 = i,\ k_2 = m-i$ we find
$$ %\be
W(s,\hat{\bf d}^{i} \cup \hat{\bf d}^{m-i}) = 
\sum_{k=0}^{s} W_{i}(k) W_{m-i}(s-k).
%\label{gauss_3}
$$ %\ee
Substituting it in (\ref{gauss_2}) and (\ref{gauss_1}) we obtain
$$
P_m^n(s) = \sum_{i=1}^m 
\sum_{k=0}^{n i - s - s_{m-i}} 
\!\!\!\!\!(-1)^{m-i}W_{i}(k) W_{m-i}(n i - s - s_{m-i} - k),
$$
where the chamber boundaries remain unchanged.
Noticing that $W_0(s)=\delta_{s,0}$ we rewrite the above expression as
\be
P_m^n(s) = W_{m}(mn-s) +
\sum_{i=1}^{m-1}
\sum_{k=0}^{n i - s - s_{m-i}} 
\!\!\!\!\!(-1)^{m-i}W_{i}(k) W_{m-i}(n i-s-s_{m-i}-k),
\label{gauss_4}
\ee
where the first term in r.h.s. contributes to {\it all} chambers.
The partial expression $P_m^n(r,s)$
for the $r$-th chamber reads
\be
P_m^n(r,s) = W_{m}(mn-s) +
\sum_{i=r}^{m-1}
\!\!\!\sum_{k=0}^{\ \ n i - s - s_{m-i}} 
\!\!\!\!\!(-1)^{m-i}W_{i}(k) W_{m-i}(n i-s-s_{m-i}-k),
\
(r-1)n \le s \le rn.
\label{gauss_5}
\ee
The relations (\ref{Gauss_known}) imply a symmetry $s \leftrightarrow mn-s$ 
leading to an 
equivalent representation
%\be
%P_m^n(s) = W_{m}(s) +
%\sum_{i=1}^{m-1}
%\sum_{k=0}^{ s - n(m-i) - s_{m-i}} 
%\!\!\!\!\!(-1)^{m-i}W_{i}(k) W_{m-i}( s - n(m-i)-s_{m-i}-k),
%\label{gauss_4a}
%\ee
%and
\be
P_m^n(r,s) = W_{m}(s) +
\sum_{i=1}^{r-1}
\!\!\!\sum_{k=0}^{\ \ s - in - s_{i}} 
\!\!\!\!\!(-1)^{i}W_{m-i}(k) W_{i}(s - ni -s_{i}-k),
\
(r-1)n \le s \le rn,
\label{gauss_5a}
\ee
and we obtain a counterpart of (\ref{gauss_4})
\be
P_m^n(s) = W_{m}(s) +
\sum_{i=1}^{m-1}
\!\!\!\sum_{k=0}^{\ \ s - in - s_{i}} 
\!\!\!\!\!(-1)^{i}W_{m-i}(k) W_{i}(s - ni -s_{i}-k).
\label{gauss_4a}
\ee
Comparison of (\ref{gauss_5}) to 
(\ref{gauss_5a}) leads to a symmetry relation
\be
P_m^n(r,s) = P_m^n(m+1-r,mn-s),
\label{symm1}
\ee
that reduces the computation of $P_m^n(s)$ to the
first $\lfloor (m+1)/2 \rfloor$ chambers.

The formulas (\ref{gauss_5a}, \ref{symm1}) 
present a closed form of the Gaussian polynomial coefficients
as a superposition of the restricted partitions discrete convolutions.
%The same time equations (\ref{gauss_5},\ref{gauss_5a}) allow to 
%reduce computation to $\lfloor (m+1)/2 \rfloor$ chambers.
Below we illustrate the above approach by deriving 
explicit expressions of $P_m^n(s)$ for small $m=3,4$  
(see also \cite{Castillo2019}). More cumbersome expressions 
required for derivation of 
$P_5^n(s)$ and $P_6^n(s)$
are presented in Appendix.

\subsection{Explicit formula for $P_3^n(s)$}
Use (\ref{gauss_4}) or (\ref{gauss_4a}) to find
$$
P_3^n(s) = W_{3}(s) +
\sum_{i=1}^{2}
\!\!\!\sum_{k=0}^{\ \ s - in - s_{i}} 
\!\!\!\!\!(-1)^{i}W_{3-i}(k) W_{i}(s - ni -s_{i}-k),
$$
or
$$
P_3^n(s) = W_{3}(3n-s) + \sum_{i=1}^2 
\sum_{k=0}^{n i - s - s_{3-i}} 
\!\!\!\!\!(-1)^{i+1} W_{i}(k) W_{3-i}(n i - s - s_{3-i} - k),
$$
%where the summands with $p \le i \le 3$ contribute to the result in the 
%chamber bounded by $s=(p-1)n$ and $s = 3n$. Thus 
and we have $W_1(s) = 1$, $s_2 = 3,\ s_1 = 1$ to obtain
\bea
P_3^n(1,s) &=& W_{3}(s) = W_{3}(3n-s) + \sum_{k=0}^{n-s-3}W_{2}(k)
-\!\!\!\!\!\sum_{k=0}^{2n-s-1} \!\!\!\!\! W_{2}(k),
\quad
0 \le s \le n,
\nonumber \\
P_3^n(2,s) &=& W_{3}(3n-s) -
\!\!\!\!\sum_{k=0}^{2n-s-1}\!\!\! W_{2}(k),
\quad
n \le s \le 2n,
\nonumber \\
P_3^n(3,s) &=& W_{3}(3n-s),
\quad
2n \le s \le 3n.
\label{P_3^n_a}
\eea
Recall that $W_2(s) = s/2+3/4 + (-1)^s/4$ and find
$$
\Sigma_2(s) = \sum_{k=0}^{s} W_{2}(k) = 
\frac{(s+1)(s+3)}{4} 
+\frac{1+\cos \pi s}{8},
$$
to arrive at the final expression 
\be
P_3^n(s) = 
\left\{
\begin{array}{lcl}
W_{3}(s),
&\quad &
0 \le s \le n,
\\
W_{3}(3n-s) - \Sigma_2(2n-s-1),
&\quad &
n \le s \le 2n,
\\
W_{3}(3n-s),
&\quad &
2n \le s \le 3n,
\end{array}
\right.
\label{P_3^n_f}
\ee
where 
\be
W_{3}(s) =  \frac{47}{72} + \frac{s}{2} + \frac{s^2}{12} + \frac{1}{8}\cos \pi s + 
\frac{2}{9}\cos\frac{2\pi s}{3}.
\label{W3}
\ee
Note that in \cite{Castillo2019} the authors 
completely describe $P_3^n(s)$ using 36 explicit expressions
while our representation requires only three formulas 
in (\ref{P_3^n_f}). It indicates that the approach based on
Sylvester-Cayley algorithm employing the idea of chambers
provides more compact and clear result.

\subsection{Explicit formula for $P_4^n(s)$}
Start with 
$$
P_4^n(s) = W_{4}(s) +
\sum_{i=1}^{3}
\!\!\!\sum_{k=0}^{\ \ s - in - s_{i}} 
\!\!\!\!\!(-1)^{i}W_{4-i}(k) W_{i}(s - ni -s_{i}-k),
$$
%$$
%P_4^n(s) = W_{4}(4n-s) + \sum_{i=1}^3 
%\sum_{k=0}^{n i - s - s_{4-i}} 
%(-1)^{i+1} W_{i}(k) W_{4-i}(n i - s - s_{4-i} - k),
%$$
and use $s_3 = 6$ to write
\be
P_4^n(s) = W_{4}(s) -
\!\!\!\sum_{k=0}^{\ \ s - n - 1} 
\!\!\!\!\!W_{3}(k)
+\!\!\!\sum_{k=0}^{\ \ s - 2n - 3} 
\!\!\!\!\!W_{2}(k) W_{2}(s - 2n-3-k)
-\!\!\!\sum_{k=0}^{\ \ s -3n - 6} 
\!\!\!\!\!W_{3}(s - 3n -6-k),
\label{P_4^n_a}
\ee
%and
%\be
%P_4^n(s) = 
%\!\!\!\!\sum_{k=0}^{n - s - 6} 
%W_{3}(n - s - 6 - k) -
%\!\!\!\!\sum_{k=0}^{2n - s - 3} 
%W_{2}(k) W_{2}(2n - s - k - 3) +
%\!\!\!\!\sum_{k=0}^{3 n - s - 1} 
%W_{3}(k)
%+ W_{4}(4n-s),
%\label{P_4^n_a}
%\ee
where
\be
W_4(s) =
\frac{2s^3+30s^2+135s+175}{288} +  
\frac{(s+5)}{32} \cos\pi s + 
\frac{1}{8}\cos \frac{\pi s}{2}+
\frac{2}{27}\left(\cos \frac{2\pi s}{3} - \cos \frac{2\pi (s+1)}{3}
\right) .
\label{W4}
\ee
%There are four chambers in this case and in the 
%first chamber ($0 \le s \le n$) all four terms in (\ref{P_4^n_a})
%are used; in the second chamber ($n \le s \le 2n$) we use last 
%three terms; in the third chamber ($2n \le s \le 3n$) last two
%terms are employed and in the fourth chamber ($3n \le s \le 4n$)
%only the last term $W_{4}(4n-s)$ survives.

\subsection{Maximal coefficient of Gaussian polynomial}
The problem of computation of the maximal coefficient of Gaussian polynomial
was addressed recently in \cite{Castillo2019} (also L. Fel, private communication).

First note that by (\ref{Gauss_known}) the maximal coefficient $p_m^n$
is given by $p_m^n = P_m^n(mn/2)$ for even $mn$ and 
$p_m^n = P_m^n((mn+1)/2)$  for odd $mn$. 
For even $m=2k$ we have to compute $P_m^n(s)$ at 
$s=kn$ that belongs to the $(k+1)$-th chamber where $kn \le s \le (k+1)n$.
For odd  $m=2k-1$ the argument of $P_m^n$ belongs to
the middle  $k$-th chamber with $(k-1)n \le s \le kn$.
Thus we have three cases
%\begin{enumerate}
%\item
%even $m=2k$, $s=kn, \quad  p_{2k}^n=P_{2k}^n(k+1,kn),$
%\item
%odd $m=2k-1$, even $n=2r,\ s=(2k-1)r,  
%\quad  p_{2k-1}^{2r}=P_{2k-1}^{2r}(k,(2k-1)r),$
%\item
%odd $m=2k-1$, odd $n=2r-1,\ s=(2k-1)(2r-1)+1, 
%\quad  p_{2k-1}^{2r-1}=P_{2k-1}^{2r-1}(k,(2k-1)(2r-1)+1).$
%\end{enumerate}
\be
p_m^n = 
\left\{
\begin{array}{lcl}
P_{2k}^n(k+1,kn), & \quad &
m=2k,%\ s=kn, 
\\
P_{2k-1}^{2r}(k,2kr-r), & \quad &
m=2k-1,\ n=2r,%\ s=2kr-r, 
\\
P_{2k-1}^{2r-1}(k,2kr-k-r+1),  & \quad &
m=2k-1,\ n=2r-1. %\ s=2kr-k-r+1.
\end{array}
\right.
\label{}
\ee
%The explicit expressions derived above for $P_3^n(s)$ and $P_4^n(s)$
%allow to find the maximal coefficients for $m=3,4$.
Use here $k=2$ and for 
$m=3$ we find $p_3^{2r} = P_{3}^{2r}(2,3r)$, and 
$p_3^{2r-1} = P_{3}^{2r-1}(2,3r-1)$; with 
$m=4$ we obtain $p_4^{n} = P_{4}^{n}(3,2n)$.

Consider first the case $m=3$. For even $n=2r$ we obtain 
$s=3r$, so that $3n-s = 3r$ and $2n-s-1 = r-1$ leading to
$P_{3}^{2r}(2,3r) = W_{3}(3r) - \Sigma_2(r-1)$.
For odd $n=2r-1$ we use $s=3r-1$, that gives 
$3n-s = 3r-2$ and $2n-s-1 = r-1$ producing
$P_{3}^{2r}(2,3r-1) = W_{3}(3r-2) - \Sigma_2(r-1)$.
After algebraic transformations we arrive at
\be
p_3^{2r} = \frac{1}{2}\left[(r+1)^2 + \cos\frac{\pi r}{2}\right],
\quad
p_3^{2r-1} = \frac{r(r+1)}{2}.
\label{p_3}
\ee
For $m=4$ we find
%the expression 
$$
P_{4}^{n}(3,2n) = W_{4}(2n) +
\!\!\sum_{k=0}^{n-1} W_{3}(k),
$$
that reduces to 
\be
p_4^{n} = \frac{(2n+5)^3}{288} +
\frac{(2n+5)}{32} +
\frac{3}{16}\cos \pi n + 
\frac{4}{27}\cos \frac{2\pi n}{3}+
\frac{4}{27} \sin \frac{\pi(8n+1)}{6}.
\label{p_4}
\ee
Turning to case $m=5$ we use $k=3$ and consider two
separate cases --  
for even $n=2r$ we obtain 
$s=5r$, so that $p_5^{2r} = P_5^{2r}(3,5r)$, while
with odd $n=2r-1$ we have $p_5^{2r-1} = P_5^{2r-1}(3,5r-2)$.
The corresponding expressions read
\bea
p_5^{2r}&=& \frac{(23r^2+69r+49)(r+1)(r+2)}{288} + \frac{425}{1728} +
\frac{3(2r+3)}{64}\cos \pi r + \frac{4}{27}\cos \frac{2\pi r}{3}
\nonumber \\
&+& 
\frac{1}{8}\left(
\cos\frac{\pi r}{2}+\sin\frac{\pi r}{2}
\right) ,
\label{p_5} \\
p_5^{2r-1}&=&
\frac{(23r^2+46r+36)r(r+2)}{288} - \frac{37}{1728} 
-\frac{1}{64}\cos \pi r -\frac{1}{8}\sin\frac{\pi r}{2}
\nonumber \\
&+&  \frac{2}{27}\cos \frac{2\pi r}{3}-
\frac{2}{27} \sin \frac{\pi(8r+1)}{6}.
\label{p_5a}
\eea
Finally, for $m=6$ we again take $k=3$ to obtain $p_6^n = P_6^n(4,3n)$ leading to
\bea
p_6^n & =& \frac{(2n+7)(66n^4+924n^3+4606n^2+9604n+12061)}{172800}
+\frac{2n^2+14n+55}{256}\cos \pi n
\nonumber \\
&+& 
\frac{1}{16}\left(
\cos\frac{\pi r}{2}-\sin\frac{\pi r}{2}
\right) + 
\frac{4}{81} \cos\frac{2\pi n}{3}
+ \frac{4}{81}\sin\frac{\pi (4n+1)}{6}
+\frac{8}{125}\left(
\cos\frac{2\pi n}{5}+\cos\frac{4\pi n}{5}
\right)
\nonumber \\
&+& 
\frac{4}{125}\left(
-\cos\frac{\pi (2n+1)}{5}+\cos\frac{\pi (4n+1)}{5}
+2\cos\frac{\pi (8n+1)}{5}
\right)
\nonumber \\
&+& 
\frac{4}{125}\left(
\sin\frac{\pi (4n+1)}{10}+2\sin\frac{\pi (8n+1)}{10}
-\sin\frac{\pi (12n+1)}{10}
\right).
\label{p_6}
\eea
It is easy to see from (\ref{p_3}-\ref{p_6}) that the
behavior of the maximal coefficient is mainly determined by the 
polynomial part of the corresponding expressions.
This observation leads us to consideration of the 
polynomial part of Gaussian polynomial coefficients.

%%%%%%%%%%%%%%%%%%%%%%%%%%%%%%%%%%%
%%%   Polynomial part
%%%%%%%%%%%%%%%%%%%%%%%%%%%%%%%%%%%%

\section{Polynomial part of Gaussian polynomial coefficients}

\subsection{Polynomial part of partition function}
It is known that a leading contribution to the partition function $W(s,{\bf d})$ 
for the set of generators ${\bf d}=\{d_1,d_2,\ldots,d_m\}$ is provided
by its polynomial part $w(s,{\bf d})$ that can be expressed through 
the Bernoulli polynomials of higher order $B_k(s,{\bf d})$ \cite{Rub04}
\be
w(s,{\bf d}) = \frac{B_{m-1}(s+\sigma_m,{\bf d})}{(m-1)!\pi_m},
\quad
\sigma_m = \sum_{i=1}^m d_i,
\quad
\pi_m = \prod_{i=1}^m d_i,
\label{Bern1}
\ee
where the Bernoulli polynomials of higher order can be defined
in umbral calculus notation as
$$
B_k(s,{\bf d}) = (s+ \sum_{i=1}^m {}^iB d_i)^k,
$$
after the expansion the replacement 
$({}^iB d_i)^k \to B_k d_i^k$ is applied, and $B_k$ denotes the $k$-th 
Bernoulli number.
This definition implies 
\be
B_k(s,{\bf d}) = \sum_{l=0}^k \binom{k}{l} s^l B_{k-l}({\bf d}),
\label{Bern_numb_expand}
\ee
and $B_k({\bf d}) \equiv B_k(0,{\bf d})$ denotes the Bernoulli number of higher order.

\subsection{Polynomial part of Gaussian polynomial coefficients}
It is instructive to consider
a polynomial part ${\cal P}_m^n(s)$ of quasipolynomial $P_m^n(s)$
\be
{\cal P}_m^n(s) =
\sum_{i=1}^{m} {\cal W}_i(m,n;s),
\label{gauss_4poly}
\ee
where ${\cal W}_i(m,n;s)$ denotes the polynomial part
of   $w_i(m,n;s)$ defined in (\ref{gauss_2}).
For $\hat{\bf d}^m$ we have $\sigma_m = s_m, \pi_m = m!$ and 
we find that 
$$
w(s,\hat{\bf d}^{i} \cup \hat{\bf d}^{m-i}) = 
 \frac{B_{m-1}(s+s_i+s_{m-i},\hat{\bf d}^{i} \cup \hat{\bf d}^{m-i})}
{(m-1)!i!(m-i)!}.
$$
Use here the replacement $s\to ni-s-s_{m-i}$ to obtain
\be
%{\cal W}_i(m,n;s) 
{\cal P}_m^n(s)=\sum_{i=1}^{m}  (-1)^{m-i}
 \frac{B_{m-1}(ni-s+s_i,\hat{\bf d}^{i} \cup \hat{\bf d}^{m-i})}
{(m-1)!i!(m-i)!}.
\label{gauss2poly1}
\ee
The polynomial part ${\cal P}_m^n(r,s)$ of the solution in $r$-th chamber reads
\be
%{\cal W}_i(m,n;s) 
{\cal P}_m^n(r,s)=\sum_{i=r}^{m}  (-1)^{m-i}
 \frac{B_{m-1}(ni-s+s_i,\hat{\bf d}^{i} \cup \hat{\bf d}^{m-i})}
{(m-1)!i!(m-i)!},
\quad
(r-1)n \le s \le rn.
\label{gauss2poly2}
\ee

%%%%%%%%%%%%%%%%%%%%%%%%%%%%%%%%%%%
%%%   Polynomial part max coefficient
%%%%%%%%%%%%%%%%%%%%%%%%%%%%%%%%%%%%

\section{Polynomial part of maximal coefficient}
The last result leads to expressions for the polynomial part of the 
maximal coefficient of the Gaussian polynomial. For small $m \le 6$ we have
\bea
&& {\cal P}_2^n(2,n) =\frac{2n+3}{4},
\nonumber \\
&&{\cal P}_3^{2r}(2,3r) =\frac{(r+1)^2}{2}+\frac{1}{36},
\quad
{\cal P}_3^{2r-1}(2,3r-1) = \frac{(3r+1)(3r+2)}{18},
\nonumber \\
&& {\cal P}_4^n(3,2n) =\frac{(2n+5)^3}{288},
\label{max_coeff_poly}\\
&&{\cal P}_5^{2r}(3,5r) =\frac{(r+1)(r+2)(23r^2+69r+49)}{288}+\frac{571}{14400},
\nonumber \\
&&
{\cal P}_5^{2r-1}(3,5r-2) = \frac{(r+2)r(23r^2+46r+36)}{288}+\frac{89}{1600},
\nonumber \\
&& {\cal P}_6^n(4,3n) =\frac{(2n+7)(66n^4+924n^3+4606n^2+9604n+7511)}{172800}.
\nonumber
\eea

\subsection{Polynomial part for even $m$}
 
For even $m=2k$ the maximal coefficient equals 
\bea
{\cal P}_{2k}^n(k+1,kn) &=&
\sum_{i=k+1}^{2k}  (-1)^{i}
 \frac{B_{2k-1}(ni-nk+s_i,\hat{\bf d}^{i} \cup \hat{\bf d}^{2k-i})}
{(2k-1)!i!(2k-i)!}
\nonumber \\
%(i=2k-j) 
&=& 
\sum_{j=0}^{k}  (-1)^{j}
 \frac{B_{2k-1}(n(k-j)+s_{2k-j},\hat{\bf d}^{j} \cup \hat{\bf d}^{2k-j})}
{(2k-1)!j!(2k-j)!}
\nonumber \\
%(\hat{\bf d}^{jk} \equiv -\hat{\bf d}^{2k-j} \cup \hat{\bf d}^{j}) 
&=& 
\sum_{j=0}^{k}  (-1)^{j}\binom{2k}{j}
\frac{B_{2k-1}(n(k-j),-\hat{\bf d}^{2k-j} \cup \hat{\bf d}^{j})}
{(2k-1)!(2k)!},
\label{gauss2poly_2k}
\eea
where we employ the relation
$B_n(s,-{\bf d}) = B_n(s+\sigma_m,{\bf d})$.
Use here (\ref{Bern_numb_expand}) to write
\bea
{\cal P}_{2k}^n(k+1,kn) = 
\frac{1}{(2k-1)!(2k)!}
\sum_{p=0}^{2k-1} \binom{2k-1}{p} n^{p}
\sum_{j=0}^{k}  (-1)^{j} \binom{2k}{j}
B_{2k-p-1}(-\hat{\bf d}^{2k-j} \cup \hat{\bf d}^{j}) (k-j)^{p}.
\label{gauss2poly_meven_sum}
\eea

\subsection{Polynomial part for odd $m=2k-1$ and even $n=2r$}
In this case the maximal coefficient reads 
\bea
{\cal P}_{2k-1}^{2r}(k,(2k-1)r) &=& 
\sum_{i=k}^{2k-1}  (-1)^{i+1}\binom{2k-1}{i}
 \frac{B_{2k-2}(2r(i-k)+r,-\hat{\bf d}^{i} \cup \hat{\bf d}^{2k-1-i})}
{(2k-2)!(2k-1)!}
\nonumber \\
%(i=2k-1-j) 
& = &
\sum_{j=0}^{k-1}  (-1)^{j}\binom{2k-1}{j}
\frac{B_{2k-2}(2r(k-1-j)+r,-\hat{\bf d}^{2k-1-j} \cup \hat{\bf d}^{j})}
{(2k-2)!(2k-1)!},
\label{gauss2poly_modd_neven}
\eea
and we find  with (\ref{Bern_numb_expand})
\bea
&&{\cal P}_{2k-1}^{2r}(k,(2k-1)r) = 
\frac{1}{(2k-2)!(2k-1)!}
\sum_{p=0}^{2k-2} \binom{2k-2}{p} r^p
\nonumber \\
&&\times
\sum_{j=0}^{k-1}  (-1)^{j} \binom{2k-1}{j}
B_{2k-p-2}(-\hat{\bf d}^{2k-1-j} \cup \hat{\bf d}^{j})(2(k-1-j)+1)^{p}.
\label{gauss2poly_modd_neven_sum}
\eea

\subsection{Polynomial part for odd $m=2k-1$ and $n=2r-1$}
We observe
\bea
{\cal P}_{2k-1}^{2r-1}(k,(2k-1)r-k+1) &=& 
\!\!\!\sum_{i=k}^{2k-1}  (-1)^{i+1}\binom{2k-1}{i}
 \frac{B_{2k-2}((2r-1)(i-k)+r-1,-\hat{\bf d}^{i} \cup \hat{\bf d}^{2k-1-i})}
{(2k-2)!(2k-1)!}
\nonumber \\
%(i=2k-1-j) 
& = &
\!\!\!\sum_{j=0}^{k-1}  (-1)^{j}\binom{2k-1}{j}
\frac{B_{2k-2}((2r-1)(k-1-j)+r-1,-\hat{\bf d}^{2k-1-j} \cup \hat{\bf d}^{j})}
{(2k-2)!(2k-1)!},
\nonumber
\\
& = &
\!\!\!\sum_{j=0}^{k-1}  (-1)^{j}\binom{2k-1}{j}
\frac{B_{2k-2}(2r(k-1-j)+r+(j-k),-\hat{\bf d}^{2k-1-j} \cup \hat{\bf d}^{j})}
{(2k-2)!(2k-1)!},
\nonumber
\eea
leading to
\bea
&&{\cal P}_{2k-1}^{2r-1}(k,(2k-1)r-k+1) =
\frac{1}{(2k-2)!(2k-1)!}
\sum_{p=0}^{2k-2} \binom{2k-2}{p} r^p
\nonumber \\
&&\times
\sum_{j=0}^{k-1} (-1)^{j} \binom{2k-1}{j}
B_{2k-p-2}(j-k,-\hat{\bf d}^{2k-1-j} \cup \hat{\bf d}^{j})(2(k-1-j)+1)^{p} .
\label{gauss2poly_modd_nodd_sum}
\eea

%\section{Bounds for maximal coefficient polynomial part}

%Consider the bounds of the maximal coefficient polynomial part.
%We observe that all terms in these expressions are positive so that
%the lower bound of them is given by the leading term in $n$. 
%This term can be found by setting $p$ equal to the upper limit in the outer sums 
%in (\ref{gauss2poly_meven_sum}-\ref{gauss2poly_modd_nodd_sum} and 
%only significantly simplified inner sum should be evaluated.
%The upper bound computation is more complex as 
%we have to estimate the double sums.

%%%%%%%%%%%%%%%%%%%%%%%%%%%%%%%%%%%
%%%  Leading term
%%%%%%%%%%%%%%%%%%%%%%%%%%%%%%%%%%%%

\section{Leading term of maximal coefficient}
Consider derivation of a general expression for the leading term
of the maximal coefficient.
The leading in $n$ term in the summand in (\ref{gauss2poly_2k})
evaluates to
$$
(-1)^{j} n^{2k-1} \binom{2k}{j}
\frac{(k-j)^{2k-1}}
{(2k-1)!(2k)!},
$$
so that the leading term $L_{n,k}$ of ${\cal P}_{2k}^n(k+1,kn)$
reads 
\be
L_{kn} = \frac{n^{2k-1}}{(2k-1)!(2k)!} 
\sum_{j=0}^{k} (-1)^{j}\binom{2k}{j} (k-j)^{2k-1},
%A_{2k-1,k},
\label{lead_term_meven}
\ee
where the sum evaluates to $A_{2k-1,k}$ being a particular case of
the Eulerian numbers of type A \cite{Carlitz1959}
\be
A_{n,k} = \sum_{j=0}^{k} (-1)^j \binom{n+1}{j}(k-j)^n.
\label{EulerianA}
\ee
%Note that the upper limit of the 
%summation can be set equal to $k-1$. 
These numbers used in a definition
of the polylogarithm function
$$
Li_{-n}(\rho) = \sum_{k=1}^{\infty} k^n \rho^k = 
\frac{1}{(1-\rho)^{n+1}}\sum_{i=0}^{n} A_{n,i} \rho^{n-i},
\quad \rho \ne 1,
$$
which in its turn relates to the 
Eulerian polynomials $H_n(\rho)$ defined through
the generating function 
$$ %\be
\frac{1-\rho}{e^t - \rho} = \sum_{n=0}^{\infty} H_n(\rho) \frac{t^n}{n!},
\quad
H_n(\rho) = (-1)^n \frac{1-\rho}{\rho} Li_{-n}(\rho).
%\label{EulerianH}
$$ %\ee
%They also satisfy the recurrence
%$A_{n,k} = k A_{n-1,k-1} + (n-k-1)A_{n-1,k}$. 
%Some sources define
%similar numbers as 
%$$
%\hat A_{n,k} = \sum_{j=0}^{k+1} (-1)^j \binom{n+1}{j}(k+1-j)^n,
%\quad
%\hat A_{n,k} = A_{n,k+1}.
%$$
It is worth to note that the function $H_n(\rho)$ can be extended into 
the Eulerian polynomials of higher order used in the 
derivation  \cite{Rub04} of the partition function $W(s,{\bf d}^m)$.
%In addition to that the Eulerian numbers  are involved in the following definitions of
%the Bernoulli numbers \cite{Jha2019} %(THEY MAY BE WRONG)
%\bea
%B_{2m} &=& \frac{(-1)^{m-1}m}{2^{2(m-1)}(2^{2m}-1)}
%\sum_{i=0}^{m-1} (-1)^{i} \sum_{l=0}^{m-i-1}
%(-1)^{l}\binom{2m}{l} (m-i-l)^{2m-1}
%\nonumber \\ & = & 
%\frac{(-1)^{m-1}m}{2^{2(m-1)}(2^{2m}-1)}
%\sum_{i=0}^{m-1}  (-1)^{i} A_{2m-1,m-i},
%\nonumber \\
%B_m &=& 
%(-1)^{m-1} \sum_{k=1}^m A_{m,m+1-k} \frac{(-1)^k}{k \binom{m+1}{k}},
%\nonumber \\
%B_{2m} %&=& \frac{(-1)^{m-1}m}{2^{2m-1}(2^{2m}-1)}
%\left[
%\sum_{l=0}^{m-1}
%(-1)^{l}\binom{2m}{l}(m-l)^{2m-1} + 
%2 \sum_{k=1}^{m-1}
%(-1)^{k}  
%\sum_{l=0}^{m-k-1}(-1)^{k+l}\binom{2m}{l} (m-k-l)^{2m-1}
%\right] \nonumber \\ 
%& = & 
%\frac{(-1)^{m-1}m}{2^{2m-1}(2^{2m}-1)}
%\left(
%A_{2m-1,m}+ 2 \sum_{k=1}^{m-1}
%(-1)^{k} A_{2m-1,m-k}
%\right).
%\label{EulerianA2}
%\eea
From (\ref{lead_term_meven}, \ref{EulerianA}) it follows that
for even $m$ the leading term $L_{mn}$ of $P_{m}^n(m/2+1,mn/2)$
reads
$$
L_{mn} = \frac{n^{m-1}A_{m-1,m/2}}{(m-1)!m!},
$$
producing $L_{2n} = n/2, \ L_{4n} = n^3/36, \ L_{6n} = 11n^5/14400$
that coincide with the leading term of the corresponding expressions
in (\ref{max_coeff_poly}).

In case of odd $m=2k-1$ we observe from (\ref{max_coeff_poly}) 
that $L_{2k-1,2r-1}=L_{2k-1,2r}$ (it is also can be seen by comparing 
(\ref{gauss2poly_modd_neven_sum}) to (\ref{gauss2poly_modd_nodd_sum})), 
so that it is sufficient to 
consider the case of even $n=2r$ only.
Using
(\ref{gauss2poly_modd_neven_sum}) %and (\ref{gauss2poly_modd_nodd_sum})
we obtain
\bea
&&L_{2k-1,2r} = 
\frac{r^{2k-2}}{(2k-2)!(2k-1)!}
\sum_{j=0}^{k-1}  (-1)^{j} \binom{2k-1}{j}
(2(k-1-j)+1)^{2k-2}.
\label{lead_modd_neven} 
%\\
%&&L_{2k-1,2r-1} = 
%\frac{1}{(2k-2)!(2k-1)!}
%\sum_{j=0}^{k-1} (-1)^{j} \binom{2k-1}{j}
%((2r-1)(k-1-j)+r-1)^{2k-2}.
\label{lead_modd_nodd}
\eea
The sum in (\ref{lead_modd_neven}) evaluates to the central MacMahon numbers
$B_{2k-1,k}$ being a particular case of the Eulerian numbers $B_{n,k}$ of type B 
% defined by the following recurrent relation for $1 \le k \le n$
that satisfy
\cite{Borowiec2015} (also see Ch.11 in \cite{Petersen2015})
$$
B_{n,k} = (2(n-k)+1)B_{n-1,k-1} + (2k+1)B_{n-1,k},
\quad
B_{n,k} = \sum_{j=0}^k(-1)^{k-j} \binom{n+1}{k-j}(2j+1)^n,
%B_{n,1} = B_{n,n} = 1,
$$
and we find 
$$
L_{2k-1,2r} = \frac{r^{2k-2}B_{2k-1,k}}{(2k-2)!(2k-1)!},
\quad \rightarrow
\quad
L_{mn} = 
\frac{(n/2)^{m-1}B_{m,(m+1)/2}}{(m-1)!m!}, 
$$
producing $L_{3,2r} = r^2/2, \ L_{5,2r} = 23r^4/288$
coinciding with the leading term of the corresponding expressions
in (\ref{max_coeff_poly}).

Collecting the above results we
obtain a general expression for the leading term of the 
maximal coefficient of the 
Gaussian polynomial through the Eulerian numbers
\be
L_{mn} = 
\frac{1}{(m-1)!m!}
\left\{
\begin{array}{lcl}
n^{m-1}A_{m-1,m/2},&
\quad &\mbox{even}\ m,
%\label{leading_term1}
\\
(n/2)^{m-1}B_{m,(m+1)/2},&
\quad &\mbox{odd}\ m, \ \mbox{even}\ n,
\\
((n+1)/2)^{m-1}B_{m,(m+1)/2},&
\quad &\mbox{odd}\ m,n.
\end{array}
\label{leading_term}
\right.
\ee

%\subsection{Upper bound of maximal coefficient}
%We found numerically the bound of the 
%Bernoulli polynomial $B_n({\bf d}^m)$ for $n>1$
%\bea
%&&|B_{n}({\bf d}^m)| < 
%\S_n({\bf d}^m) n B_n, 
%\quad \mbox{even} n,
%\nonumber \\
%&&|B_{n}({\bf d}^m)| < 
%\S_n({\bf d}^m) n B_{n-1}, 
%\quad \mbox{even} odd,
%\nonumber \\
%\eea
%where $\S_n({\bf d}^m) = \sum d_i^n$.

%%%%%%%%%%%%%%%%%%%%%%%%%%%%%%%%%%%
%%%   Acknowledgments
%%%%%%%%%%%%%%%%%%%%%%%%%%%%%%%%%%%%

\vskip0.5cm
\section*{Acknowledgments}
The author thanks L.Fel for
drawing attention to the Gaussian polynomials problem and numerous fruitful discussions.

%%%%%%%%%%%%%%%%%%%%%%%%%%%%%%%%%%%
%%%   Appendices
%%%%%%%%%%%%%%%%%%%%%%%%%%%%%%%%%%%%

%%%%%%%%%%%%%%%%%%%%%%%%%%%%%%%%%%%
%%%   Appendix A
%%%%%%%%%%%%%%%%%%%%%%%%%%%%%%%%%%%%
\appendix
\section{Computation of $P_5^n(s)$ and $P_6^n(s)$\label{appendix2}}
\renewcommand{\theequation}{\thesection\arabic{equation}}
\setcounter{equation}{0}

As it follows from (\ref{gauss_1}) and (\ref{gauss_2}) the solution 
for $m=5$ can be constructed
using only three basic expressions. Namely, we have
\bea
W(s,\hat{\bf d}^2 \cup \hat{\bf d}^3)  &=&
\frac{s^4+18s^3+112s^2+279s+220}{288}+\frac{37}{1728}+
\frac{2s+9}{64}(-1)^s + \frac{2}{27} \cos\frac{2\pi s}{3},
\nonumber \\
W(s,\hat{\bf d}^1 \cup \hat{\bf d}^4)  &=&
\frac{s^4+22s^3+166s^2+495s+465}{576}+\frac{25}{3456}+
\frac{2s+11}{128}(-1)^s
\nonumber \\
&+& \frac{2}{27} \sin\frac{(4s+1)\pi }{6}
+ \frac{1}{16}\left(
\cos\frac{\pi s}{2}+\sin\frac{\pi s}{2}
\right),
\nonumber \\
W(s,\hat{\bf d}^5) &=&
\frac{s^4+30s^3+310s^2+1275s+1687}{2880}+\frac{41}{86400}+
\frac{2s+15}{128}(-1)^s + \frac{2}{27} \cos\frac{2\pi s}{3}
\nonumber \\
&+& 
+ \frac{1}{16}\left(
\cos\frac{\pi s}{2}+\sin\frac{\pi s}{2}
\right) +
\frac{2}{25}\left(
\cos\frac{2\pi s}{5}+\cos\frac{4\pi s}{5}
\right).
\label{A1}
\eea 
Then we use (\ref{gauss_2}) to find $w_i(5,n;s)$ to insert into (\ref{gauss_1}) or 
(\ref{gauss_1r})
for explicit expression of $P_5^n(s)$ or $P_5^n(r,s)$ respectively.

The solution for $m=6$ is produced
based on the four following expressions
\bea
W(s,\hat{\bf d}^3 \cup \hat{\bf d}^3)  &=&
\frac{6 s^5+180 s^4+2020 s^3+10440 s^2+24299 s+19650}{25920} 
+ \frac{s+6}{64}(-1)^s 
\nonumber \\
&+& \frac{6s+34}{243} \cos\frac{2\pi s}{3}
+ \frac{4}{243}\sin\frac{\pi (4s+1)}{6},
\nonumber \\
W(s,\hat{\bf d}^2 \cup \hat{\bf d}^4)  &=&
\frac{12 s^5+390 s^4+4720 s^3+26130 s^2+64458 s+54275}{69120}+
\frac{2s^2+26s+75}{512}(-1)^s 
\nonumber \\
&+& \frac{2}{81} \cos\frac{2\pi s}{3}
+ \frac{2}{81}\sin\frac{\pi (4s+1)}{6} + 
\frac{1}{32}\left(
\cos\frac{\pi s}{2}+\sin\frac{\pi s}{2}
\right),
\nonumber \\
W(s,\hat{\bf d}^1 \cup \hat{\bf d}^5)  &=&
\frac{6 s^5+240 s^4+3560 s^3+24000 s^2+71325 s+70888}{86400}+
\frac{s+8}{128}(-1)^s
\nonumber \\
&+& \frac{2}{81} \cos\frac{2\pi s}{3}
+ \frac{2}{81}\sin\frac{\pi (4s+1)}{6}
+ \frac{1}{16}\sin\frac{\pi s}{2} +
\frac{2}{125}\left(
\cos\frac{2\pi s}{5}+\cos\frac{4\pi s}{5}
\right)
\nonumber \\
&+&
\frac{2}{125}\left(
\cos\frac{\pi (2s+1)}{5}+3\cos\frac{\pi (6s+1)}{5}
+2\cos\frac{\pi (8s+1)}{5}
\right)
\nonumber \\
&-&\frac{2}{125}\left(
2\sin\frac{\pi (8s+1)}{10}+\sin\frac{\pi (12s+1)}{10}
+3\sin\frac{\pi (16s+1)}{10}
\right),
\nonumber \\
W(s,\hat{\bf d}^6) &=&
\frac{12 s^5+630 s^4+12320 s^3+110250 s^2+439810 s+598731}{1036800}+
\frac{6s^2+126s+581}{4608}(-1)^s 
\nonumber \\
&+& 
\frac{6s+61}{486} \cos\frac{2\pi s}{3}
+ \frac{2}{243}\sin\frac{\pi (4s+1)}{6}
+\frac{1}{18} \cos\frac{\pi s}{3}
\nonumber \\
&+& 
\frac{1}{32}\left(
\cos\frac{\pi s}{2}+\sin\frac{\pi s}{2}
\right) +
\frac{2}{125}\left(
\cos\frac{2\pi s}{5}+\cos\frac{4\pi s}{5}
\right)
\nonumber \\
&+&
\frac{2}{125}\left(
\cos\frac{\pi (2s+1)}{5}+3\cos\frac{\pi (6s+1)}{5}
+2\cos\frac{\pi (8s+1)}{5}
\right)
\nonumber \\
&-&\frac{2}{125}\left(
2\sin\frac{\pi (8s+1)}{10}+\sin\frac{\pi (12s+1)}{10}
+3\sin\frac{\pi (16s+1)}{10}
\right).
\label{B1}
\eea 
Now employ (\ref{gauss_2}) to obtain $w_i(6,n;s)$ and then
use it in (\ref{gauss_1}) or (\ref{gauss_1r})
for explicit expression of $P_6^n(s)$ or $P_6^n(r,s)$ respectively.

%%%%%%%%%%%%%%%%%%%%%%%%%%%%%%%%%%%%%%%%%%%%%%%%%%%%%%%%
%%%%%%%%%%%%%%%%%%%%%%%%%%%%%%%%%%%%%%%%%%%%%%%%%%%%%%%%
%%%%%%%%%%%%%%%%%%%%%%%%%%%%%%%%%%%%%%%%%%%%%%%%%%%%%%%%

%%%%%%%%%%%%%%%%%%%%%%%%%%%%%%%%%%%
%%%   References
%%%%%%%%%%%%%%%%%%%%%%%%%%%%%%%%%%%%

\end{document}